\documentclass[11pt]{elsarticle}

\usepackage{comment}
\usepackage{geometry}
\usepackage{natbib}
\usepackage{soul}
\setcitestyle{round,authoryear}
\usepackage{hyperref}
\usepackage[english]{babel}
\addto\extrasenglish{
    
}
\usepackage{array}
\usepackage{multirow}
\usepackage{tabularx}
\usepackage{bm}
\usepackage{amssymb}
\usepackage{booktabs}
\usepackage{caption}
\usepackage{longtable}
\usepackage{amsfonts,amsmath}
\usepackage{xcolor}
\usepackage[linesnumbered]{algorithm2e}
\allowdisplaybreaks
\usepackage{bbm}
\usepackage{tikz}
\usepackage{float}
\usepackage{subfig}
\usepackage{amsmath}
\usepackage[normalem]{ulem}
\usetikzlibrary{shapes.misc, positioning}
\usetikzlibrary{arrows.meta}
\usepackage{mdframed}
\usepackage{lscape}

\def\bx{{\boldsymbol x}}  
\newcounter{cdcounter}  
\definecolor{Gray}{gray}{0.85}
\newcolumntype{a}{>{\columncolor{Gray}}c}

\begin{document}

\begin{frontmatter}



\title{Synthesizing mixed-integer linear programming models from natural language descriptions}




\author[label1,label2]{Qingyang Li}
\ead{ql5@student.unimelb.edu.au}
\author[label1,label2]{Lele Zhang}
\ead{lele.zhang@unimelb.edu.au}
\author[label3]{Vicky Mak-Hau}
\ead{vicky.mak@deakin.edu.au}

\affiliation[label1]{organization={The University of Melbourne},
              city={Melbourne},
              postcode={3010},
             state={Victoria},
             country={Australia}}
\affiliation[label2]{organization={ARC Training Centre in Optimisation Technologies, Integrated Methodologies, and Applications (OPTIMA)},
             country={Australia}}
\affiliation[label3]{organization={Deakin University},
              city={Burwood},
              postcode={3125},
             state={Victoria},
             country={Australia}}

\begin{abstract}

Numerous real-life problems arising in various domains such as production planning, resource allocation, scheduling, and transportation planning and management, can be effectively formulated and solved using Mixed-Integer Linear Programming (MILP) models. However, the transformation of real-world decision-making problems into MILP models heavily relies on expertise in operations research and mathematical optimization, which restricts non-experts' accessibility to MILP. 
To address this challenge, we propose a framework for automatically formulating MILP models from unstructured natural language descriptions of decision problems, which integrates Large Language Models (LLMs) and mathematical modeling techniques.
This framework consists of three phases: i) identification of decision variables, ii) classification of objective and constraints, and iii) finally, generation of MILP models.

In this study, we propose a knowledge representation structure that includes relational knowledge, inheritable knowledge, and so on. Specifically, we present a constraint classification scheme and a set of constraint templates that can guide the LLMs in synthesizing a complete MILP model from a given natural language description of a decision-making problem. After \textit{fine-tuning} LLMs based on numerous examples of constraint and objective descriptions, our approach can, from the natural language description, identify and synthesize logic constraints in addition to classic demand and resource constraints. The logic constraints have not been studied in existing work. 


To evaluate the performance of the proposed framework, we extend the NL4Opt dataset with more problem descriptions and constraint types, in particular logic constraints and binary variables, and with the newly created dataset, we compare our method with one-step model generation methods offered by LLMs like ChatGPT and Google Bard. The experimental results reveal that with respect to the accuracies of generating the correct model, objective, and constraints, our method which integrates knowledge representation with LLMs significantly outperforms the others.
Taking ChatGPT for example, our framework accurately formulates 86.67\% problems and 99.32\% constraints compared to ChatGPT's 26.67\% and 82.99\% respectively. 
The three-phase framework that we developed is a prototype system for synthesizing MILP models from natural language descriptions and has a great potential to capture more constraints for more complex MILPs. Our framework opens up opportunities for developing training tools for operations research practitioners and has the potential to be a powerful tool for automatic decision problem modeling and solving in practice.

\end{abstract}



\end{frontmatter}

\section{Introduction}

In the last few decades, mathematicians and computer scientists have collectively improved the way business decisions are made, through understanding the requirements of the end-users and creating computerized systems to optimize their decisions. 
Many real-life decision-making problems are combinatorial optimization problems that can typically be formulated as Integer Programming (IP)-family of problems (such as Pure Binary Integer Programming (BIP), Mixed-Integer Linear Programming (MILP) problems). 
CO problems are about making the best decisions given a large number of (and sometimes exponentially many) options. For example, to make decisions on ``yes/no'' questions such as whether to build a facility at a location, to allocate a machine to a job, or if a customer is to be visited immediately before another customer on a trip. Or, to make decisions on a quantity, for example, how many tables to make, how many acres of land to be used to grow apples, or how much money to invest in a stock. These decisions are often made with the objective to maximize profit or minimize cost, subject to business requirements or needs, such as limitations due to resource availability, quantities to be balanced, demands to be supplied, quality to be conformed to, or business logic to be observed. 

Optimizing business decisions through formulating mathematical programming models and developing solution methodologies to solve them has brought huge social and economic benefits. 
For instance, Ford Motor developed an Integer Programming (IP) model to shorten the planning process, which optimized global procedures and saved USD 250M \citep{Chen2009}. Baosteel applied IP for production planning, and just their main Shanghai plant alone, reduced annual carbon monoxide emissions by over 500,000 tons \citep{Edelman2013}. As of 2020, (M)ILP was the technique used in one out of five Franz Edelman Award Finalists \citep{Gorman2020}. In fact, businesses of all sizes can benefit from optimized decision making. 
 
 In the last six decades or so, the focus of research amongst the Operation Research community has been on the theories, solution methodologies, and applications of IP-family of problems, see, e.g., \citep{50years}. Countless MILPs have been formulated for different combinatorial optimization problems in the past, and no doubt many more to come. An MILP comprises decision variables (binary, integer, and/or continuous), an objective function (a linear combination of parameters and decision variables), and constraints of different classes. Even though an MILP can have exponentially many constraints, the number of types of MILP constraints is limited. With recent advances in Natural Language Processing (NLP), in particular Large Language Models (LLMs), the idea of automating the formulation of MILPs has inspired much new research. 

 Our research is motivated by the desire to find out how LLMs can be utilized to  ``translate'' a natural language description of a combinatorial optimization problem into a mathematical model, in particular, through constraint recognition and model synthesis.

\subsection{Literature review}

Surprisingly, research on the automatic synthesis of mathematical programming models from natural language descriptions has only become active in recent years. 
There has been a lot of research thus far on methods that synthesise mathematical programming models from data (observations). They made use of different techniques that allow the synthesis of different types of Mathematical Programming (MP) models, such as utilizing MILP formulation for the linear and nonlinear constraint synthesis problems \citep{pawlak17automatic},
local search for MILP model synthesis \citep{sroka18oneclass},
inductive logic programming for integer nonlinear programming \citep{kumar19acquiring},
Grammatical Evolution for MILP constraint synthesis \citep{pawlak21grammatical}.
These data-driven constraint learning methods do not focus on the context of the problem represented by the model, and hence, they are more applicable to the constraints that are difficult to formulate directly. However, one concern is that the interpretability of constraints learned from the data may be poor and not easily understood by humans \citep{fajemisin23optimization}. Another limitation of these approaches is the availability of data, i.e., whether a sufficient number of solution examples can be collected for good performance. Notice that the only guarantee of a correct MILP being ``learned'' from data is when all ``positive samples'' (feasible solutions) are present in the system. However, in most real-life combinatorial optimization problems, even finding one feasible solution is hard. If all feasible solutions are known, then one can just evaluate the objective function value for all the solutions to find the optimal solution. 

Research on automatically solving natural language mathematical word problems has been conducted since the 1960s. An intelligent system capable of solving mathematical word problems should be able to process the natural language description of the problem, allow different users to describe the problem in different ways, and minimize control over the user's input \citep{mukherjee08reviewa}. Many approaches have been proposed for transforming natural language descriptions of mathematical word problems into sets of equations, which are mostly based on semantic parsing, text similarity, predefined templates for sets of equations, or deep learning \citep{zhang19gap}. 

The rise of LLMs opens up more possibilities for equation synthesis of mathematical word problems. LLMs can perform new tasks with only instruction and no training examples (zero-shot learning), or with a few training examples (few-shot learning) \citep{kojima23large}.
\citet{zong23solving} evaluated the performance of Generative Pre-trained Transformer 3 (GPT-3, \citeauthor{brown20language}, \citeyear{brown20language}) in three types of tests:  classifying mathematical word problems that can be expressed as two linear equations; extracting equations from word problems; and generating word problems. The accuracy of extracting equations increases with the number of examples provided in the prompt, and the accuracy of the fine-tuned GPT-3 is higher than that of the few-shot learning and zero-shot learning. GPT-3 with zero-shot learning performs well on classifying word problems for most classes, and their experiment results provide some evidence that supplying examples of the same class in the prompt for few-shot learning for extracting equations can be helpful. However, in contrast to our approach, they do not incorporate classification as a prior step to equation extraction in their framework. Notice that problem classification is not always helpful in equation extraction. Compared to our approach, these authors fine-tuned GPT-3 used for equation extraction rather than for the classification task with a limited number of valid responses. 

In contrast to their study, we classify constraints rather than problems since constraint classification is more helpful for extracting MILP models.
In integer programming for combinatorial optimization (CO) problems, sometimes different problems can share the same constraints. For example, many of the scheduling problems have similar requirements to routing problems. Sometimes, a CO problem can have characteristics of two or more well-known problem classes. With different combinations of constraints, there are infinitely many possibilities of CO problems that can be formulated as MILPs, hence classifying problems may not be a practical solution. On the other hand, the number of constraint types is limited, and therefore classifying the type of a constraint is much easier.



Recent technological advancements within the field of natural language processing have opened up more possibilities on the topic of automatic synthesis of mathematical programming models from natural language descriptions.
Extracting mathematical models from text descriptions can be regarded as a semantic parsing task, i.e., extracting machine-interpretable meaning representations from a natural language text \citep{ramamonjison23nl4opt}. 
\citet{OFOGHI2023110980} presented an  ``ontology'' for MILP, however some elements are unclear. For example, coefficients were said to be  parametric, but there is no mention that the decision variables are parametric as well. The paper also states that ``order is of importance for indices that represent parametric coefficients, so the correct sequence is preserved''. As MILP practitioners, our understanding is that the order of appearance of the set of decision variables and their associated coefficients is usually arbitrary. Further, if the purpose of the proposed ontology is for automatic model generation, then without seeing the framework, it is hard to ascertain if the ontology is useful. If the purpose is only to demonstrate a possible machine representation of MILPs, then this work may not be new, as there are so many open source and commercial MILP solvers, each of them may have developed an ontology at some point. 


One of the earliest works in this area is by \cite{islam21automatic}, incorporating both the fields of natural language processing and linear programming. Given a structured problem description, their method first classifies the objective as maximizing or minimizing using Support Vector Machine (SVM). It then tokenizes sentences using the Python NLTK package and extracts the numbers using NLTK POS-tagger. However, there are obvious limitations to their method, namely the structured problem description and the requirement that all sentences included must be simple. Information about model components (e.g., number of variables, number of constraints) and parameter values (e.g., coefficients of variables) must appear in the textual description in a specified order. Moreover, the method can only handle constraints that represent upper bounds.

With the popularity of LLMs growing, a number of recently developed automated modeling frameworks have utilized LLMs in the methodology architecture for converting natural language descriptions into mathematical models. 
To formulate and solve the MILP problem from a structured natural language description, \citet{ahmaditeshnizi23optimus} developed an LLM-based agent called OptiMUS. 
The method uses LLM to build the mathematical formulation in one step by using a prompt containing user input that contains additional ``hints'' far beyond the natural language description of the problem.
OptiMUS is also able to exploit the LLM to transform the mathematical model into solver code, debug it, and check the validity of the solution. The study shows the potential to automate multiple stages of solving optimization problems by integrating LLMs and solvers. 
However, it requires users to use as input a representation called SNOP (Structured Natural Language Optimization Problem) that explicitly specifies the decision variables, the objective function, and all the parameters used to compose the objective and constraints. Whilst this greatly reduces the challenge of automatic modeling, writing an accurate input is difficult for users without mathematical programming-related expertise. 
Typically, a domain expert does not describe his/her decision problem in such an ideally structured language style.
The fact is that different people are likely to describe the same problem in different ways. In contrast to their work, our paper is concerned with automatic formulation from unstructured text descriptions of MILP problems. 

In 2022, the NL4Opt competition \citep{ramamonjison23nl4opt} for extracting linear programming formulations from natural language has attracted wider attention to this research topic and has led to significant progress in this field. 
\citet{ramamonjison22augmenting} proposed an automatic formulation method, OptGen, that is capable of handling unstructured LP problem descriptions with different types of constraints, and developed and tested it on a dataset containing 1101 problem instances. 
In the first stage (entity identification), the fine-tuned transformer acts as an entity tagger that tags the keywords denoting variables, objective function, and constraints from the problem description. 
In the second stage, a set of intermediate representation declarations of the objective and constraints are generated from the recognized entities, creating a bridge from contextual descriptions to context-independent mathematical formulations. An intermediate representation declaration is a sequence of tokens derived from the problem description to preserve the contextual information of the problem, indicating the objective or constraint (and its type) in a structured way. This transformation is based on predefined templates for different types of objectives and constraints. 
Finally, the transformation of intermediate representation declarations to mathematical expressions is performed by a parser. The authors created a modeling platform that allows user interaction using the OptGen architecture described above as the underlying framework. 
In addition, they developed the first dataset of linear programming word problems (NL4Opt) that comprises instances in a number of application domains.
Since then, a number of research groups utilize ensemble learning to improve the first-stage outcomes \citep{he22linear, doan22vtccnlp, wang23opd, ning23novel}. Data pre-processing and augmentation are used to improve the second-stage outcomes \citep{gangwar23highlighting, jang22tag, ning23novel}. \citet{prasath23synthesis} used GPT-3 to generate new instances by making parameter changes and context changes to the original problems in order to increase the size of the training dataset. The authors reported that in most cases, the new problem descriptions generated were not semantically correct. Despite this, the models were trained on the augmented data but evaluated using only the original NL4Opt dataset, and it reportedly led to improved performance. 

These studies are the first steps toward the goal of automatic linear programming modeling from natural language descriptions. The average number of decision variables in the NL4Opt instances is 2.08 and the average number of constraints is 2.83. The constraints are confined to a number of basic constraint types, and neither logic constraints nor binary variables are accounted for. Whilst the problem descriptions covered a wide range of domains, the difficulty of the mathematical models is limited to simple linear programming (even though the decision variables in most problem instances should in fact be integers). These are typically simple LP problems taught in the first week of a tertiary optimization subject, or even in some secondary school mathematics textbooks. 


To transform any combination of Boolean logic operations (such as \texttt{AND}, \texttt{OR} and \texttt{Implication}) into MILP constraints, \citet{wang23synthesizing} developed a syntax-guided synthesis method based on Domain Specification Language (DSL). The input specification is a logical condition consisting of propositions representing linear expressions connected by Boolean logic operators, e.g., $(X_8 = 0 \rightarrow X_{11} = 0) \land (X_8 = 0 \rightarrow X_{17} = 0)$, which means \texttt{if $X_8=0$, then $X_{11}$ and $X_{17}$ must be $0$}. This type of specification can be seen as an intermediate representation between natural language descriptions to mathematical expressions. 

Advances in LLMs make it possible to develop a new workflow for automated modeling.
Inspired by the NL4Opt competition, \citet{tsouros23holy} developed a framework called Holy Grail 2.0, consisting of similar steps, i.e., extracting the semantic entities of the optimization problem (as proposed in \citep{dakle23ner4opt}), identifying the relationships between the entities, formulating the problem, and generating runnable code, but the framework heavily utilizes LLMs in all steps involved. To assess the ability of the method to formulate problems, they presented four levels of abstraction for natural language problem descriptions, (from all model components being stated explicitly to the closest way that someone not familiar with optimization modeling would describe a problem.) 
For developing a diverse instance set of MILP word problems in the future, the abstraction level of a problem description can be considered as a characteristic of the instance, providing guidance for creating new instances that are more challenging in terms of translating natural language descriptions into mathematical models.

\citet{li23large} developed another LLM-based framework called OptiGuide for several application scenarios in supply chain optimization, with textual queries as input and interpretations of optimization solutions as output. In this system, LLMs are used to transform the user's questions described in natural language (e.g., can we use café S1 for roastery D1?) into new constraints added to the MIP model (e.g., \texttt{m.addConstr(x["S1", "D1"]==1)}) or updates to the parameters in the MIP model in code form.
In this approach, a prompt consists of simple instructions and examples, i.e., pairs of questions and code answers, in addition to the question asked by the user. The approach requires developing, for each scenario, a repository of question-code pairs, some of which form a set for the selection of examples in the prompt, and the others serve as a test set. The performance of few-shot learning is highly affected by the appropriateness of the example selection in the prompt. In contrast, our approach uses classification plus fine-tuning to avoid the shortcomings of few-shot learning and to shorten the length of the input to LLMs by avoiding providing a large number of examples in the prompt.

In all studies discussed above, LLMs are employed as translators in the transformation from natural language to mathematical models. However, the OPRO (optimization by prompting) method developed by \citet{yang23largea} uses LLMs as optimizers to improve optimization solutions. They describe optimization problems in natural language but do not formally formulate them. Instead of using solvers to obtain solutions, the method guides the LLM to iteratively generate new solutions by providing previously generated solutions and their evaluations as examples in the prompt, since LLMs are capable of identifying patterns from examples. This method is only suitable for small-scale optimization problems.

In contrast to the aim of transforming a given problem description into a mathematical model, some research focuses on the generation of models that satisfy the specifications required by users. It provides us with another approach for building instance sets of MILP word problems, namely, generating problem descriptions from diverse models. \citet{almonacid23automatic} tentatively demonstrated the ability of GPT-3.5 to generate optimization models (represented in the Minizinc modeling language) that satisfy the required features of the model as specified by the user via prompts, such as whether the domains of the variables are open or bounded, and whether constraints are imposed.

\subsection{Contributions and outline of the paper}


The main contributions of this paper include:
\begin{itemize}
    \item proposing a framework to automatically build MILP models from unstructured natural language descriptions of optimization problems, which integrates LLMs and mathematical modeling techniques.
    This framework consists of three stages as shown in \autoref{fig:flow1}: i) identification of decision variables, ii) classification of the objective and constraints, and iii) constraint generation and supplementation. 
    In contrast to previous work, this approach is able to handle binary variables and logic constraints. It has a great potential to capture more constraints for more complex MILPs.
    This framework has been shown to outperform LLMs when they are asked to directly generate the MILP models, with the former generating correct models in 50\% to 60\% more problem instances.
    Our framework provides opportunities for the development of training tools for operations research practitioners. It also has the potential to become a powerful tool for automatic decision problem modeling and solving in practice, making optimal decision making more accessible to non-experts.
    \item developing an instance set of MILP word problems, comprising 30 problem descriptions and their MILP models. In contrast to the first LP word problem instance set, NL4Opt, our instance set involves binary variables and includes four types of logic constraints.
    \item proposing a constraint recognition method by fine-tuning LLMs as constraint classifiers. We extend the constraint descriptions in the NL4Opt dataset (the dev set) by adding descriptions for logic constraints, annotating their classes, and fine-tuning LLMs on this modified dataset.
    \item creating templates for constraints to guide LLMs in generating the mathematical formulation of the MILPs.
\end{itemize}

The rest of the paper is organized as follows. \autoref{sec:constraint_type} describes the types of constraints we currently consider and compares them with the constraint classification of the NL4Opt competition. \autoref{sec:methodology} illustrates the methodology of the framework for automatically converting natural language descriptions into mathematical models. \autoref{sec:numerical} presents the numerical experiment, including dataset creation and results comparing the performance of two variants of our framework and two chatbots in automatically synthesizing models. Finally, \autoref{sec: conclusion} summarizes the paper and proposes some future work.







\section{Constraint Classification} \label{sec:constraint_type}


\begin{table}[!ht]
    \small
    \begin{tabular}{ p{0.7cm}p{6.8cm}p{2.2cm}p{4.5cm}p{0.7cm}}
    \hline
    \textbf{Type No.}&   \textbf{Constraint type defined in this work}& \textbf{Mathematical inequality} &\textbf{Constraint type defined in NL4Opt competition} &\textbf{Type No.}\\
    \hline
    1 & Upper bound on single variable & $ x_i \leq b $ &Upper bound &2\\
    2 & Upper bound on sum of variables & $\sum_{i} x_i \leq b $ &Sum constraint &1\\
    3 & Upper bound on weighted sum of variables & $\sum_{i} a_i x_i \leq b $ &Linear constraint &4\\
    4 & Upper bound on proportion & $ x_j \leq c \sum_{i} x_i $ &Ratio control constraint & 5 \\ 
    5 & Lower bound on single variable & $ x_i \geq b $ &Lower bound &3\\
    6 & Lower bound on sum of variables & $\sum_{i} x_i \geq b $  &Sum constraint &1\\ 
    7 & Lower bound on weighted sum of variables & $\sum_{i} a_i x_i \geq b $ &Linear constraint &4\\
    8 & Lower bound on proportion & $ x_j \geq c \sum_{i} x_i $ &Ratio control constraint &5\\
    9 & Comparison constraints & $ d x_i \leq x_j $ &Balance constraint Type-1 (if $d \neq 1$) or Type-2 (if $d =1$) &6/7\\
    10 & If A then B/ if not B then not A/ B if A & $y_A\leq y_B$ & & \\
    11 & Exactly one of A and B/ either A or B but not both (or neither) & $y_A + y_B = 1$ & & \\
    12 & At least one of A and B/ if not A then B/ either A or B or both & $y_A + y_B \geq 1$ & & \\
    13 & At most one of A and B/ if A then not B/ either A or B or neither (but not both) & $y_A + y_B \leq 1$ & & \\
    \hline
    \end{tabular}
    \caption{Constraint types defined in this work and those in the NL4Opt competition. $x_i$ and $x_j$, $i, j \in \{1,2, \ldots, n\}$, denote decision variables, $a_i$, $b$ and $d$ are non-negative constants, and $c$ is a constant $\in (0, 1]$. We use a binary decision variable $y \in \{0,1\}$ to represent the true/false value of a statement (e.g., if an event will take place), with, e.g., $y_A = 1$ indicating Statement A is true (Event A does take place), and $y_A = 0$ otherwise. }\label{tab:constraint class}
\end{table}

Even though the types of constraints are limited, it is non-trivial to consider all constraint types for automatic modeling. This work, as a proof-of-concept, focuses on the categories of constraints: \textit{basic constraints} included in the NL4Opt dataset and some \textit{logic constraints}. 

The NL4Opt competition classified the constraints in their problem instances into 7 types, e.g., upper bound, ratio control and sum constraints, and here we refine them into 9 types. See \autoref{tab:constraint class}. 
We deconstructed the 99 problem descriptions from the NL4Opt (development) dataset, and obtained 292 descriptions of constraints. These constraints can be covered by the types 1 to 9 in \autoref{tab:constraint class}. 

However, the problems in the NL4Opt dataset only involve continuous or integer variables and the nine types of constraints, and no binary variables or logic constraints are included. 
In this study, we consider 4 types of logic constraints, each involving two binary variables. 
To extend the NL4Opt dataset, we create descriptions of logic constraints belonging to the types 10 to 13 in \autoref{tab:constraint class} and problem instances involving these types of constraints.
Due to limited existing datasets for MILP modeling, this study relies on the NL4Opt dataset, which does not contain equality constraints. Therefore this study does not consider equality constraints except for the logic constraint type 11 that we have added.

In what follows, we simplify our notation and use $ax$ to represent  $\boldsymbol{a} \cdot \bx$, where $\boldsymbol{a}$ and $\bx$ are vectors of coefficients and decision variables respectively. We also use $\sum_j$ to represent $\sum_{\forall j}$. 

\subsection{Constraints for resource, demand, proportions and comparison} \label{ssec:resource_demand} 

Constraints that capture resource restrictions (static upper limit) and minimum demand requirements (static lower limit) are commonly used in MILP. They are expressed in the form of $a x \leq b$ (e.g., types 1-3 in \autoref{tab:constraint class}) and $a x \geq b$ (e.g., types 5-7) respectively. If the requirement is to equate amounts, we should use $ax =b$. However, this equality constraint type is not addressed in the NL4Opt dataset, and hence, this study focuses on the inequality constraints. 

In some cases, upper or lower bounds on proportions are required, e.g., types 4 and 8 in \autoref{tab:constraint class}. For example, the proportion of product $j$ must not be more than $c$ of the total (mathematically, $x_j \leq c \sum_j x_j$ with $0 < c \leq 1$). Or, {the sum of a product set ${\cal N}_1$ must be at least $c$ of the sum of another set ${\cal N}_2$ ($\sum_{j\in{\cal N}_1} x_j \geq c\sum_{j\in{\cal N}_2} x_j$).}


A comparison between two quantities can be given as $ d x_i \leq  x_j$ for $d \in \mathbb{R}_{+} $. 
If $d>1$, it implies that the amount of $j$ must be at least $d$ times the amount of $i$ (e.g., $d=3, 3 x_i \leq  x_j$), or that the amount of $i$ cannot exceed $1/d$ of the amount of $j$ (e.g., $d=3, x_i \leq (1/3) x_j$). 
If $d<1$, it implies that the amount of $j$ must be at least $d$ of the amount of $i$ (e.g., $d=1/3, (1/3) x_i \leq  x_j$), or that the amount of $i$ cannot exceed $1/d$ times the amount of $j$ (e.g., $d=1/3, x_i \leq 3 x_j$).
In the NL4Opt competition dataset, the constraint descriptions of type 9 involve only the first and third expression styles. Therefore, we do not need to specifically define another type of comparison constraints represented by $ x_i \leq d x_j$. Nevertheless, in terms of using LLMs to formulate comparison constraints, there is no impact since LLMs have the flexibility to use the data mentioned in the constraint descriptions as the coefficients of variables.

\subsection{Set packing/partitioning/covering constraints} \label{ssec:set}

Set packing/partitioning/covering constraints are sub-types of resource/conservation/demand constraints respectively. These constraints require at most/exactly/at least one of $n$ options to be selected. Constraints of type 11/12/13 in \autoref{tab:constraint class}) are special cases of set partitioning/covering/packing constraints. These options are represented by a set of $n$ binary variables $\bm y \in \{ 0, 1 \}^n$. Then set packing/partitioning/covering constraints are denoted as $\bm a\cdot \bm y \{ \leq, =, \geq \} 1$ respectively for $\bm a \in \{ 0, 1 \}^n$. 
If the right-hand side (RHS) of the constraint is an integer constant $b > 1$, it is the weighted set packing/partitioning/covering constraint. 
Set packing and set covering constraints as well as weighted versions of them are also included in type 2 ($\sum_{i} x_i \leq b $) and type 6 ($\sum_{i} x_i \geq b $) in \autoref{tab:constraint class} if the variables $x_i$'s are all binary. 


\subsection{Logic constraints} \label{ssec:logic}

 In real-life applications, sometimes just capturing quantities (using continuous or general integer variables) is not enough, we need to capture logic requirements as well.
Binary variables are commonly used to convert logical relations in a combinatorial optimization problem into linear constraints, known as logic constraints. Consider an example from the dataset we created; \texttt{if Haus Toys makes trucks, then they will not make trains}, can be represented as the nonlinear constraint $\texttt{trucks}*\texttt{trains}=0$, where the integer variables $\texttt{trucks}, \texttt{trains}\in\mathbb{Z}_{\geq0}$ refer to the number of trucks and trains they plan to produce, respectively. 
To capture the same business requirement as a linear programming model, one can introduce binary variables $\texttt{bi\_trucks}, \texttt{bi\_trains}\in\{0,1\}$ indicating whether Haus Toys makes trucks and trains respectively. The binary variable takes the value of 1 if they make trucks/trains, and 0 otherwise. Then the logic constraint can be represented as the linear inequality $\texttt{bi\_trucks} + \texttt{bi\_trains} \leq 1$.

In general, one can use a binary variable $y_A\in \{0,1\}$ to capture if $x_A=0$ or $x_A>0$ for a continuous or integer variable $x_A$. To establish the link between the binary variable and its corresponding continuous or integer variable, we impose two \textit{linking constraints} as follows.
\begin{align}
    x_A \leq M \times y_A, \label{eq:link_con1 in sec2}\\
    y_A \leq M \times x_A. \label{eq:link_con2 in sec2}
\end{align}
Inequality \eqref{eq:link_con1 in sec2}, for a sufficiently large value $M$, enforces $y_A=1$ if $x_A>0$. Inequality \eqref{eq:link_con2 in sec2} ensures $y_A=0$ if $x_A=0$ (here, $M$ is needed only if $x_A$ may take values less than 1). 



Consider two statements, A and B, with truth values represented by binary variables $y_A, y_B \in \{0,1\}$ respectively; for example, $y_A=1$ means statement A is true and 0 otherwise. This study considers the following logic conditions.
\begin{itemize}
    \item The \textit{If-then} condition: ``If A then B" (type 10), equivalent to ``If not B then not A", is given by $y_A\leq y_B$. 
    Since ``not A" can be represented by $1-y_A$, ``If not A then B" (type 12), equivalent to ``At least one of A and B", is satisfied by $1-y_A\leq y_B$.
    ``If A then not B" (type 13), equivalent to ``At most one of A and B", is satisfied by $1-y_B\geq y_A$.
    \item The \textit{Exclusive or} condition (type 11): ``Either A or B (but not both)" is given by $y_A + y_B=1$, which is equivalent to ``Exactly one of A and B".
    \item The \textit{Inclusive or} condition (type 12): ``Either A or B or both" is given by $y_A + y_B \geq 1$, which is equivalent to ``At least one of A and B".
\end{itemize}

\section{Methodology}
\label{sec:methodology}

\begin{figure}[!ht]
  \centering
  \includegraphics[width=1\linewidth]{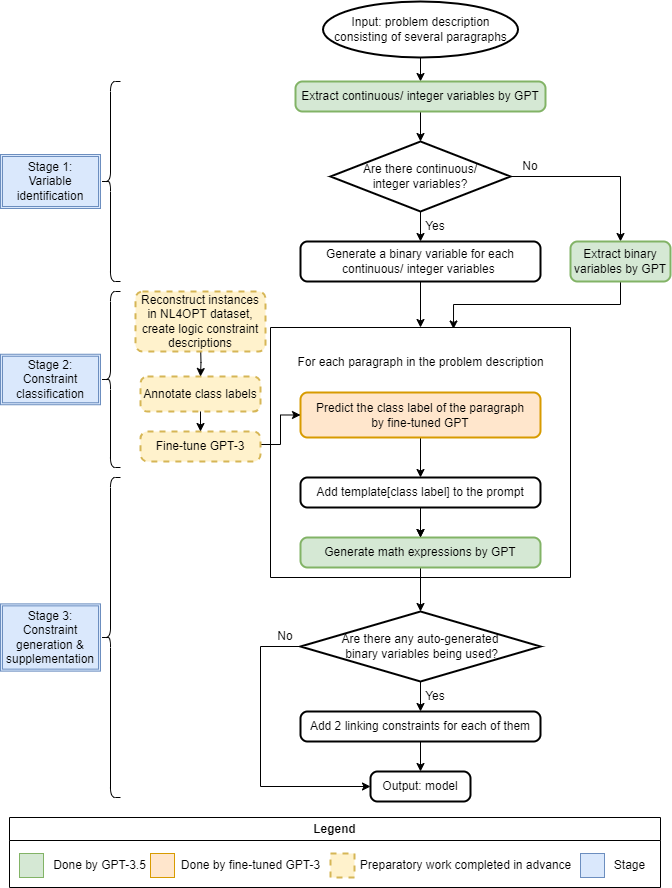}
  \caption{Main steps of the automatic modeling approach with GPT}
  \label{fig:flow1}
\end{figure}
\begin{landscape}
\begin{figure}[!ht]
  \centering
  \includegraphics[width=1\linewidth]{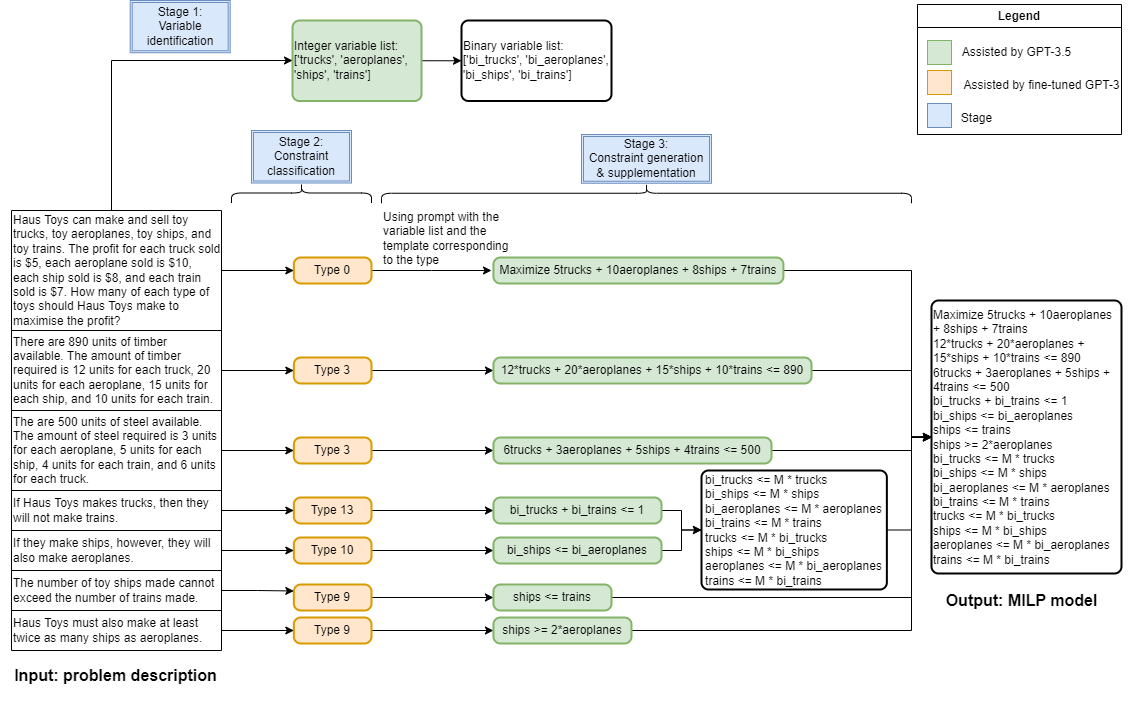}
  \caption{Implementing the automatic modeling approach on an example from our dataset}
  \label{fig:flow2}
\end{figure}
\end{landscape}

The automatic modeling approach proposed in this paper consists of three stages and two preparation modules.
The three stages are 1) variable identification, 2) constraint classification, and 3) constraint generation and supplementation. See \autoref{fig:flow1}.
We assume that a full problem description is provided as input by the user and it is presented in the form of paragraphs, each representing an objective or a constraint. See an example in \autoref{fig:flow2}.
Given the problem description, the approach first determines all decision variables, then for each paragraph, classifies the objective function or constraint and generates a mathematical formula using the corresponding template, and finally for logic constraint(s) if any, supplements linking constraints, to the final model as the output.

The first preparation module fine-tunes an LLM to be the classifier of Stage 2, and the second preparation module, empowered by the knowledge representation, creates a set of constraint templates for generating mathematical constraints in Stage 3. Both modules are completed before executing the automatic modeling procedure.
The main flow of the proposed approach is depicted in \autoref{fig:flow1}, whereas an example of the problem description input, the modeling process, and the model output is shown in \autoref{fig:flow2}.

\subsection{Stage 1: Variable identification}
\label{subsec:stage1}
 
The first stage of the modeling approach is to determine the decision variables and is composed of two steps. The first step is to identify all continuous and/or integer variables from the full problem description by an LLM. If the problem does not involve any types of these variables, the LLM is asked to recognize all binary variables. 
In that case, the problem is considered to involve only binary variables and will be formulated as a pure binary linear programming model.
The answer of the LLM should consist only of the variable names.
We instruct the LLM that a continuous or integer variable needs to be called with the name of the quantity that it represents. For example, the integer variable representing \texttt{the quantity of vanilla cakes to be made by a bakery} is named \texttt{vanilla\_cakes}. The name of a binary variable must start with \texttt{bi\_} and be linked to the name of the entity it represents, for example, \texttt{bi\_vanilla\_cakes}. The benefit of this naming convention is that it facilitates comparison of the model formulated with the true model, eliminating the effects of the order in which the variables appear. 


The second step of Stage 1 is to generate a corresponding binary indicator variable for each of the continuous or integer variables; if its value is greater than 0 then the binary variable is 1, otherwise 0. This step is the basis for constructing the logic constraints later in the automatic modeling approach. It is important to note that the binary variables generated in this step are associated with the corresponding continuous or integer variables, and so in Stage 3, for those of them that are used in the formulation, we need to add linking constraints to the model in order to associate each of them with the corresponding continuous or integer variable.
Note that all variables in the dataset used in this paper are non-negative.

\subsection{Stage 2: Constraint classification}
\label{subsec:stage2}

For each paragraph of the given problem description, which is an objective function description or a constraint description, the second stage classifies the constraint type from the type set in \autoref{tab:constraint class} for each description using a fine-tuned LLM. 

LLMs are pre-trained on large amounts of data, which allows them to respond to a given prompt containing only a few examples or no examples. 
However, for the constraint classification in this framework, we allow LLMs to rely on more than just prompts and to learn more examples. Thus, we fine-tune LLMs as the classifier of Stage 2 to recognize constraints, so that there is no need to provide examples in the prompts. 
Fine-tuning trains an already trained model by using a specific dataset to adapt it to a customized application scenario. It can improve the accuracy and response speed of LLMs to specific tasks. 
More details about the dataset and hyperparameters used for fine-tuning can be found in \autoref{sec:numerical}. Once the dataset and hyperparameters are set up, the LLM fine-tuning, such as for GPT models and PaLM, can often be done automatically.

For the classification by LLMs, if the name of each class (i.e., type) contains more than one token, OpenAI suggests that it may affect the output quality of the trained model. 
Therefore, in the dataset for fine-tuning, we use the type numbers rather than type names in \autoref{tab:constraint class} to annotate the objective function and the 13 constraint types in each prompt-completion pair. 
This allows us to use a single objective description or a single constraint description as the prompt to a fine-tuned LLM, and the fine-tuned LLM responds with its corresponding type number as the completion.

\subsection{Stage 3: Constraint generation and supplementation}
\label{subsec:stage3}

\subsubsection{Constraint generation}
\label{subsubsec:stage3generation}

The third stage is to generate mathematical formulas for the objective function and constraints under the guidance of the templates. 
In Stage 2, we have identified whether a paragraph in the given problem description described an objective function or a type of constraint. 
For each paragraph, Stage 3 uses the predefined template corresponding to the type label classified in Stage 2 as part of the prompt to guide the LLM for \textit{template matching}, and further generates the mathematical expression using the variables defined in Stage 1. Here, template matching refers to matching the variables identified to the variable placeholders in the template and then substituting the identified variables into the placeholders for variables in the formula templates. In order to standardize the form of the answers generated by LLMs, we also use some instructions in combination with the template to form the prompt. Examples of the prompts for generating the objective function and some types of constraints can be found in \ref{app:prompt}. 

A template for a constraint type explains the meaning of such constraint and the form of the corresponding mathematical formula. For logic constraint types, the templates also include the usual expression of the constraint in natural language to facilitate correct template matching by LLMs. For example, given binary variables $y_A$ and $y_B$ each taking the value of 1 are denoted as statements A and B each being true, the logic constraint belonging to type 10 is described in natural language as ``If A (is true) then B (is true)", which corresponds to the mathematical expression $y_A \leq y_B$. Given a description, for example, \texttt{if Product A is produced then Product B is produced}, the template for type 10 constraints allows LLMs to learn how to substitute the variables into the formula template. In this example, if the expression uses the ``$\leq$" sign, then the binary variable representing the state of the proposition following ``if" (whether Product A is produced) should be placed to the left-hand-side of the inequality sign, and the binary variable indicating the state of the proposition following ``then" (whether Product B is produced) should be placed to the right-hand-side of the inequality sign. 

The main contents of the templates are shown in \autoref{tab:templates}. Since LLMs are better at handling text prompts, we tend to use natural language in templates to describe the meaning and form of the formula, rather than using complex mathematical symbols directly as prompts.

\begin{table}[!h]
    \caption{Templates for the objective function and 13 types of constraints. Consider Statements A and B with truth value represented binary variables $y_A$ and $y_B$ respectively - with 1 representing the statement is true and 0 otherwise. } \label{tab:templates}
    \begin{tabular}{r|p{9.2cm} | p{6.0cm}}
    \hline
    \textbf{No.}&\textbf{Meaning} & \textbf{Formula template}\\
    \hline
    0 & The objective function should specify a direction of optimization, either to maximize or minimize. Therefore, the answer should only include the words ``Maximize" or ``Minimize" and a linear expression.  & sum of terms consisting of a variable multiplied by one or more constant coefficients\\
    1 & This constraint represents an upper bound on a single decision variable. & variable $<=$ constant\\
    2 & This constraint represents an upper bound on the sum of decision variables. & sum of variables $<=$ constant \\
    3 & This constraint represents an upper bound on the weighted sum of decision variables.  & sum of variables multiplied by their weight $<=$ constant\\
    4 & This constraint indicates that the upper bound on a single variable is a constant proportion of the sum of all variables. & one variable $<=$ proportion * sum of all variables\\
    5 & This constraint represents a lower bound on a single decision variable. & variable $>=$ constant\\
    6 & This constraint represents a lower bound on the sum of decision variables.  & sum of variables $>=$ constant \\
    7 & This constraint represents a lower bound on the weighted sum of decision variables.  &  sum of variables multiplied by their weight $>=$ constant \\
    8 & This constraint indicates that the lower bound on a single variable is a proportion of the sum of all variables.  & one variable $>=$ proportion * sum of all variables\\
    9 & This constraint is a comparison constraint between two variables. & c * x $<=$ y, where x and y are variables and c is a positive constant.\\
    10 & In natural language descriptions, this type of constraint often contains a format like ``If A then B" or ``If not B then not A". & $y_A <= y_B$\\ 
    11 & In natural language descriptions, this type of constraint often contains a format like ``(Exactly) one of A and B" or ``Either A or B (but not both)". & $y_A + y_B = 1$\\
    12 & In natural language descriptions, this type of constraint often contains a format like ``At least one of A and B" or ``If not A then B" or ``Either A or B or both". & $y_A + y_B >= 1$\\
    13 & In natural language descriptions, this type of constraint often contains a format like ``At most one of A and B" or ``If A then not B" or ``Either A or B or neither (but not both)". & $y_A + y_B <= 1$\\
    \hline
    \end{tabular}
\end{table}

\subsubsection{Constraint supplementation}
\label{subsubsec:stage3supplementation}

After each constraint or objective function in the problem description has been formulated, if the problem is not a pure binary problem (one that involves binary variables only), the final task of Stage 3 is to introduce two additional constraints for each binary indicator variable used in the generated mathematical model. 
Specifically, if any of the binary variables automatically generated in the second step of Stage 1 are used in the formulation, then the model needs to add two constraints for each used binary variable to link it to the corresponding continuous or integer variable. Even if the following two types of linking constraints are not explicitly mentioned in the problem description, they capture some business rules. 

For example, \texttt{if an investor decides not to invest in a particular stock, then the amount invested in this stock cannot be greater than 0}. 
To ensure that if the binary variable is equal to 0, then its corresponding continuous or integer variable should be at most 0, we need to add the following constraint to the model,
\begin{equation}
    \text{variable\_name} \leq M \times \text{bi\_variable\_name}, \label{eq:link_con1}
\end{equation} 
where $M$ is a constant representing a sufficiently large value. The binary variable is required to be equal to 1 if its corresponding continuous or integer variable is greater than 0. 
However, if this continuous or integer variable has an upper bound, $M$ will be replaced by this upper bound, and the upper bound constraint on this variable generated earlier in Stage 3 can be replaced by inequality \eqref{eq:link_con1}.

In addition, the binary variable should be 0 if its corresponding continuous or integer variable is equal to 0. In other words, if the binary variable is equal to 1 then its corresponding continuous or integer variable should be greater than 0. To enforce this, the following constraint should be added to the model,
\begin{equation}
    \text{bi\_variable\_name} \leq M \times \text{variable\_name},
    \label{eq:link_con2}
\end{equation} 
where again, $M$ is a sufficiently large constant.
Taking an example in the investment context, the above constraint implies that \texttt{if an investor decides to invest in a particular stock, then this means that he needs to invest at least some amount that is greater than 0}.

The objective function and constraints translated from the problem description plus automatically supplemented linking constraints (if needed) constitute the final MILP model as the output shown in \autoref{fig:flow2}.

\section{Numerical Experiment}
\label{sec:numerical}


The proposed multi-stage framework requires the assistance of LLMs for variable identification in the first stage, constraint classification (with fine-tuned LLMs) in the second stage, and constraint generation in the third stage. In this numerical experiment, we employ two well-known LLMs: GPT and PaLM 2. The first variant, which we refer to as the multi-stage approach with GPT, fine-tunes GPT-3 to classify constraints and uses GPT-3.5 to identify variables and generate constraints. The second variant, which we refer to as the multi-stage approach with PaLM, adopts PaLM 2 in all the stages.  

We compare two variants of our approach with two baseline models that utilize zero-shot prompting to directly instruct the following chatbots in directly converting problem descriptions to mathematical models.
\begin{itemize}
    \item ChatGPT:  AI-powered chatbot, powered by GPT-3.5 and developed by OpenAI.
    \item Bard: chat-based AI tool, powered by PaLM 2 and developed by Google.
\end{itemize}


\subsection{Dataset and configuration for fine-tuning}
\label{subsec:Dataset for fine-tuning}



In our experiments, we fine-tuned two LLMs, GPT-3 and PaLM 2, to implement the classification prediction for a given segment description of an objective or a constraint. They were later used in Stage 2 when testing the entire framework. 
The first dataset on which we conduct the numerical experiment is for fine-tuning LLMs, namely training and validating classifiers. 
It consists of instances of objective descriptions and constraint descriptions as well as their annotated type labels. 

Each instance of NL4Opt dataset has a paragraph describing an LP problem. We deconstructed the 99 problem descriptions from the NL4Opt (development) dataset, and obtained 391 descriptions of objectives and constraints. Each problem description has only one objective function. These 292 constraints can be covered by nine types, i.e., types 1 to 9 in \autoref{tab:constraint class} and \autoref{tab:constraint number}. However, the problems in the original NL4Opt dataset do not involve logic constraints. To extend the decomposed NL4Opt dataset, we create a total of 183 descriptions of logic constraints involving two binary variables belonging to four types, i.e. types 10 to 13 in \autoref{tab:constraint class} and \autoref{tab:constraint number}. The numbers of descriptions contained in these types are shown in \autoref{tab:constraint number}.

The dataset used for fine-tuning thus contains 574 pairs of descriptions and manually annotated type numbers, which are then divided into the training and validation sets automatically by the data preparation tool of OpenAI. 
Due to the relatively small size of the current dataset, the number of descriptions for some classes is small. There are only 7 descriptions for type 6, and they were all split into the training set, resulting in missing classes in the validation set which is not ideal. We therefore manually moved one description of type 6 from the training set to the validation set.
As a result, the training set contains 464 descriptions, and the validation set contains 110 descriptions.
For the purpose of a fair comparison, we use the same training set and validation set to fine-tune PaLM.

\begin{table}[!ht]
    \centering
    \caption{The dataset used in fine-tuning LLMs for constraint classification}\label{tab:constraint number}
    \small
    \begin{tabular}{ p{1cm}|p{7cm}|p{3.6cm}|p{3.6cm}}
    \hline
    \textbf{Type No.}&   \textbf{Type name}&\textbf{No. of descriptions from NL4Opt}& \textbf{No. of descriptions we create for tuning}\\
    \hline
    0& Objective function&99 &-\\
    1 & Upper bound on single variable &20 &-\\
    2 & Upper bound on sum of variables & 12&-\\
    3 & Upper bound on weighted sum of variables & 93&-\\
    4 & Upper bound on proportion & 8&-\\ 
    5 & Lower bound on single variable & 36&-\\
    6 & Lower bound on sum of variables & 7&-\\ 
    7 & Lower bound on weighted sum of variables & 59&-\\
    8 & Lower bound on proportion & 14&-\\
    9 & Comparison constraints & 43&-\\
    10 & If A then B/ if not B then not A/ B if A & -&46\\
    11 & Exactly one of A and B/ either A or B but not both (or neither) & -&40\\
    12 & At least one of A and B/ if not A then B/ either A or B or both & -&49\\
    13 & At most one of A and B/ if A then not B/ either A or B or neither (but not both) & -&48\\
    \hline
    \end{tabular}
\end{table}


After preparing the dataset, we need to select the base model and the hyperparameters to start the automatic fine-tuning.
GPT-3 series has four base models, which were the only models of OpenAI that could be fine-tuned at the time we conducted our experiments. For proof-of-concept purposes, we fine-tuned the base model called \textit{ada} for this specific text classification task, which is suitable for simple classification and is the fastest and least expensive.
We also tuned the base model \textit{text-bison-001}, which is the only base model of the PaLM API that can be tuned so far.
Hyperparameters used for fine-tuning GPT-3 and PaLM, such as the number of epochs and batch size, can be found in \autoref{tab:hyperparameters}. Finding better classifiers using different configurations of the hyperparameters, namely hyperparameter tuning, can be left as further research for future work, which is not the focus of this paper for proof-of-concept purposes.

\subsection{Test set for evaluating the framework} 
\label{ssec:evaluation dataset}


The second dataset on which we conduct numerical experiments is for evaluating the whole framework of automatic modeling.
In contrast to NL4Opt dataset, we develop a new test set, comprising 30 MILP word problems to introduce the four types of logic constraints mentioned in \autoref{sec:constraint_type}. 
Each problem consists of its natural language description and a mathematical model (for evaluating the accuracy). For each problem, the model has a single objective function and several constraints. In the problem description, each paragraph describes either an objective function or a constraint. The dataset includes a total of 30 objective descriptions and 147 constraint descriptions. 
In terms of both natural language contexts and mathematical formulations, the problems within the new dataset showcase diversity. 
The descriptions of these problems span different domains, including production planning, resource allocation, and transport. These problems include different combinations of more constraint types. 

Readers may find examples of the prompts given to ChatGPT and Bard as well as those used in the multi-stage framework in \ref{app:prompt}. 
All data used in this paper is available as a GitHub repository (\url{https://github.com/yangyangyang777/Datasets-for-automatic-modeling.git}).

\subsection{Evaluation metrics} \label{ssec:metrics}

We evaluate the performance of the framework using two instance-level and expression-level formulation accuracy metrics, ACC1 and ACC3. For evaluating the classification performance of Stage 2, we also consider another expression-level classification accuracy, ACC2. The formulas for the three metrics are given below.
\begin{align}
    \text{Model generation: } \text{ACC1}& = \frac{\sum_{i=1}^N T_i}{N},\label{eq:ACC1}\\
    \text{Classification of objectives and constraints: } \text{ACC2}& = \frac{\sum_{i=1}^N TC_i}{\sum_{i=1}^N D_i},\label{eq:ACC2}\\
    \text{Generation of objectives and constraints: } \text{ACC3}& = \frac{\sum_{i=1}^N TE_i}{\sum_{i=1}^N D_i}.\label{eq:ACC3}
\end{align}
Here, the number of problem instances in the test set is denoted by $N$. The total number of the objective descriptions and constraint descriptions contained in the problem description of instance $i$ is denoted as $D_i$. 
For a given instance $i$, $T_i$ takes the value of $1$ if the generated model is completely accurate (including linking constraints if any) and $0$ otherwise. 
For the expression-level metrics, $TC_i$ is the number of objective and constraint descriptions in instance $i$ that are correctly classified, and $TE_i$ is the number of objective and constraint descriptions that are correctly formulated. 

We do not specifically investigate the accuracy of variable identification in this section. If the decision variables of instance $i$ are incorrectly identified (variables are missing and/or unneeded variables are defined), {it would be reflected by $T_i$. Specifically, if some decision variable(s) of instance $i$ is missing, then the model is incorrect, at least partially wrong, and hence $T_i = 0$. If some redundant variable(s) is defined, the model is incorrect and $T_i=0$ unless the redundant variable does not appear in any of the objective and constraint expressions but only in the variable list.
}

\subsection{Fine-tuning results} \label{ssec:Fine_tuning_results}

We use classification accuracy to initially evaluate the performance of the fine-tuned GPT-3 and PaLM on the validation set, before integrating them with the framework and evaluating the whole framework on the test set.
Similar to the ACC2 metric, classification accuracy here refers to the number of objective and constraint descriptions that are correctly classified divided by the total number of descriptions in the fine-tuning validation set.

The accuracy of the fine-tuned GPT-3 reached 99.09\%. Only one instance in the validation set containing 110 instances was misclassified. The description CD1 was classified as class 4 (Upper bound on proportion), but the actual class is class 2 (Upper bound on the sum of decision variables). In the training set, both types 2 and 4 have less than 10 constraints. The small number of samples available for learning may have led to this misclassification. OpenAI hypothesized that the model would perform better with at least 100 examples per class in theory. Expanding the dataset for constraint classification should be a future task, especially for those classes with few instances, so that LLMs can learn patterns from more instances.  

\begin{quote}
Constraint description 1 (CD1) \texttt{Bold Tycoon decides to invest his money in GICs and index ETF. Bold Tycoon wants to invest $\$10,000$ in total.}
\end{quote}

\begin{quote}
Constraint description 2 (CD2) \texttt{A waste treatment company must remove waste using a large container or a medium container. There must be at most 65 total containers.}
\end{quote}


The tuned PaLM obtained an accuracy of 98.18\%. Two instances in the validation set were misclassified. But by analyzing the results it appears to be more reliable than the fine-tuned GPT-3 on constraint classification.
The description CD1 was classified as type 3 (Upper bound on weighted sum of variables), which makes sense somehow since the actual class (upper bound on sum of decision variables) is a special case of type 3. The misclassification of CD2 was the same as that of CD1. 
Type 3 contains far more constraints than any other constraint type since it is very common in practice, and this imbalanced fine-tuning dataset may have an impact on strict constraint classification. However, fuzzy (not uniquely correct) classifications similar to the above is acceptable in terms of modeling.

\subsection{Modeling results} \label{ssec:Results}

We assessed the two variants of the multi-stage framework on the new instance set consisting of 30 problems and compared them with two baselines. The experimental results in \autoref{tab:resultstable} reveal that the multi-stage approaches which integrate knowledge representation with the two LLMs significantly outperform the one-step model-generation methods by ChatGPT and Bard with respect to the accuracies at both the model and expression levels. Between the two LLMs, the two methods powered by the GPT models outperform those powered by PaLM. 

\begin{table}[!h]
    \caption{Results of the multi-stage framework (utilizing GPT models or PaLM) proposed in this paper and two baseline models, ChatGPT and Bard.} \label{tab:resultstable}
    \begin{tabular}{p{4.99cm}|p{2.9cm} |p{2.7cm} | p{2cm}| p{2cm}}
    \hline
    &\textbf{Multi-stage approach with GPT} &\textbf{Multi-stage approach with PaLM} & \textbf{ChatGPT} 
    & \textbf{Bard}\\
    \hline
    ACC1 & 0.8667 & 0.7 & 0.2667 & 0.1667 \\
    ACC2 & 0.9944 & 0.9887 & - & -\\
    ACC3 & 0.9774 & 0.9435 & 0.8475 & 0.7627 \\
    No. of incorrect models & 4 &9 & 22 & 25\\
    No. of incorrect objectives & 3 & 4 &  2 & 3 \\
    No. of incorrect constraints & 1 & 6 & 25 & 40 \\
    No. of models missing linking constraints & 1 & 0 & 17 & 17 \\
    \hline
    \end{tabular}
\end{table}

\subsubsection*{Cases missed by the multiple-stage approach}

The multi-stage approach using GPT did not formulate the correct models for four problem instances in the test dataset. Among them, three models were incorrect because the coefficients of the variables in the generated objective functions were wrong. Taking one of the three models as an example, the objective is to maximize the total profit. The coefficients for computing the profits should be expressed as the differences between selling prices and costs, yet the generated objective function used selling prices directly as coefficients.
The multi-stage approach with PaLM and the one-step approach with Bard made the same mistake. Surprisingly, ChatGPT correctly formulated the objective function for the profit maximization. 
The errors in the other two objective functions were due to the failure to perform unit conversions on the coefficients. 
Taking one of these two problems as an example, the variables are in minutes. The data referring to the variable coefficients in the objective function description are in minutes per kilometer. They cannot be used as coefficients directly since their inverse (in kilometers per minute) should be used.
The multi-stage approach with PaLM, ChatGPT and Bard all failed to notice the inconsistency in the units and derived incorrect models for the two problems. In the future, to improve the multi-stage approach, we may consider including an additional calculation phase for unit conversion.

Our framework with GPT misclassified one logic constraint as a non-logic constraint, giving this classification model an accuracy ACC2 of 99.44\% on the test dataset. As a result, the constraint was incorrectly formulated. We took a closer look at the problem instance, which is to determine the number of orders from three manufacturers. The description of the misclassified constraint is given below.
\begin{quote}
Constraint description 3 (CD3): \texttt{If the store decides to order chairs from manufacturer A, they must also order at least 10 chairs from manufacturer B.}
\end{quote}

CD3 belongs to the logic constraint of the type ``if A then B” (i.e., type 10 in \autoref{tab:constraint class}). However, the second half of the sentence is close to the description of the basic lower bound constraint (i.e., type {5} in \autoref{tab:constraint class}). 
The chance of accurately classifying it would be improved if the method makes use of the information from other constraint descriptions of this problem. From the whole problem description, it can be noticed that one type of decisions to be made is to determine the number of orders from different manufacturers (integer variables). Specifically, another constraint description is ``\texttt{Each order from manufacturer A will include 15 chairs, while each order from B and C will include 10 chairs.}”. With that, CD3 is equivalent to ``\texttt{If the store decides to order chairs from manufacturer A, they must also order chairs from manufacturer B.}”. Classifying and modeling this equivalent description is less challenging. This example suggests that it requires global information as well as logical reasoning for modeling challenging constraint descriptions. 
Similarly, our framework with PaLM misclassified CD3 as another type of logic constraint thus failing to formulate it correctly. Google Bard and ChatGPT also failed to model it correctly.

Since the multi-stage approach with GPT formulated CD3 as a non-logic constraint, its mathematical expression did not use binary variables. As a result, the linking constraint attached to the binary variables for the logic constraint was missed in the final model of the problem, which explains the missing linking constraint in \autoref{tab:resultstable} (second column, last row). 
On the other hand, the approach with PaLM classified CD3 as a logic constraint using two binary variables therefore the corresponding linking constraints were added. 
For this problem, ChatGPT and Bard did not correctly formulate any of the logic constraints and did not add the linking constraints.

In addition to the above mentioned, more mistakes were made by the multi-stage approach using PaLM and these were mainly due to missing coefficients of variables, and incorrect direction of the inequality (ignoring the specification of upper or lower bounds in the template).

\subsubsection*{Cases missed by ChatGPT and Bard}

ChatGPT's mistakes are mainly in the formulation of constraints; in particular, failing to formulate 20 logic constraints, often due to its failures in identifying the binary variables. These constraints were formulated as logic constraints in terms of integer and continuous variables. But even for those problems where it correctly generated expressions for logic constraints, ChatGPT missed the linking constraints, leading to incorrect final models. For other types of constraints, common mistakes include having the wrong coefficients for the variables and missing variables.

Furthermore, ChatGPT missed three constraints in three problem descriptions. For another problem instance, ChatGPT generated a constraint that is not mentioned in the given problem description. Both of the errors were impossible to occur within our framework.

In comparison, Bard's performance was worse than that of ChatGPT. Bard only generated correct models for five of the thirty instances. Bard was more prone to miss constraints than ChatGPT, missing 7 constraints in total.
Bard failed more times than ChatGPT on the formulation of both the basic constraints and the logic constraints. ChatGPT made 4 mistakes for basic constraints, and Bard made 14 due to using wrong variables, coefficients, the direction of inequality in the constraints, etc. 
There is also a special class of errors that Bard made, using unexpected variables outside of the given problem description to formulate a constraint. For one problem instance with seven continuous variables, Bard correctly identified and defined the seven continuous variables, but used an eighth continuous variable in a constraint whilst this variable was not defined by Bard and not included in the problem description.

Bard was less reliable in formulating logic constraints. One important finding is that Bard could only correctly formulate logic constraints for pure binary programming models. For MILP problems, it did not once identify binary variables thus failing to correctly formulate logic constraints, not to mention the linking constraints. In contrast, ChatGPT successfully identified binary variables for a few MILP problems.  

\subsubsection*{Variable identification}

For every instance in the test dataset, the variable identification procedure of our framework successfully defined all the variables needed to build its mathematical model, especially the binary variables, which facilitates high accuracy in the subsequent formulation phase. LLMs by themselves are often unaware of binary variables involved in the problem description and have less knowledge about logic constraints. 


\section{Conclusion}
\label{sec: conclusion}


In this paper, we propose a three-stage framework to automatically synthesize MILP models from unstructured natural language descriptions of decision problems.
This framework can employ LLMs to perform tasks in all stages, namely the identification of decision variables, the classification of objective and constraints, and the generation of the MILP model. 
We propose a constraint recognition method for Stage 2 by using the modified NL4Opt dataset to fine-tune LLMs as constraint classifiers. 
We build constraint templates to instruct LLMs to generate formulas for Stage 3.

The framework can be adapted to make it capable of handling more types of constraints for more complex MILP problems. For the purpose of proof-of-concept, we first extended the application of the framework to four types of logic constraints, in contrast to previous studies that can only handle classic demand and resource constraints. 
To fill this gap, we developed a dataset of the MILP word problems with logic constraints. 
The results of testing the framework on this dataset show that our method which integrates knowledge representation significantly outperforms one that uses LLMs alone. The accuracy of generating the correct model can be improved by about 50\% to 60\%.

Some directions can be explored as future work. 
For the framework presented in this paper, we assume that a full problem description has been provided as input. In the next step, we will adjust the framework to another scenario, assuming that we elicit from the user information sufficient to construct an optimization model through a question-and-answer dialogue system. In this context, the two modules of eliciting the problem description and mathematical formulation will not be separated, but may be more interactive. Such a problem elicitation approach may facilitate automatic partitioning between constraint descriptions, which is done manually in this paper. 
This paper focuses on the translation of natural language problem descriptions to mathematical models, and further translation to modeling languages for integration with solvers is future work. 
The performance of the LLMs integrated in this framework can potentially be improved by adding appropriate few-shot examples to the prompts. 
The dataset needs to be expanded to include not only more contexts but also more constraint types and optimization problem types. Testing the capability of LLMs to generate new MILP word problems can facilitate exploring systematic generation of instances.

\section*{Acknowledgment}

This work was partially funded by the Australian Research Council, Australia through the Discovery Project 2022 (grant number DP220101925). It was partially funded by the Australian Government through the Australian Research Council Industrial Transformation Training Centre in Optimisation Technologies, Integrated Methodologies, and Applications (OPTIMA), Project ID IC200100009.


\bibliography{sample_base}

\newpage

\appendix

\section{Appendix}

\subsection{Hyperparameters for fine-tuning}
\label{app:hyper}

Hyperparameters used for fine-tuning GPT-3 and PaLM are shown in \autoref{tab:hyperparameters}. \cite{openai} suggests that the default values of the hyperparameters they have picked work well in various use cases. The recommended batch size is about 0.2\% of the number of instances in the training set, which we set to 1.

\begin{table}[!h]
    \centering
    \caption{Hyperparameters used for fine-tuning GPT-3 and PaLM} \label{tab:hyperparameters}
    \begin{tabular}{p{4.2cm}|p{2cm}|p{4.2cm}|p{2cm} }
    \hline
    \textbf{Hyperparameter for fine-tuning GPT-3} &\textbf{Value} &\textbf{Hyperparameter for fine-tuning PaLM} &\textbf{Value}\\
    \hline
    Epochs & 4 (default) & Epochs & 20 \\
    Batch size & 1 & Batch size & 24 \\
    Learning rate multiplier & default  & Learning rate & 0.02  \\
    \hline
    \end{tabular}
\end{table}

\subsection{Prompt examples} \label{app:prompt}

\sethlcolor{lightgray}


Taking the problem in \autoref{fig:flow2} as an example, a few prompt examples are depicted in \autoref{fig:prompt to ChatGPT and Bard} - \autoref{fig:prompt for type9}. Problem-specific texts have been \hl{shaded}. The prompt examples for constraint generation (\autoref{fig:prompt for type3} - \autoref{fig:prompt for type9}) are displayed in the order in which these constraints appear in the problem description.

\begin{table}[h!]
    \centering
        \caption{An example of the prompts to ChatGPT and Bard}
    \label{fig:prompt to ChatGPT and Bard}
    \begin{mdframed}
        \small
        \setlength{\parskip}{5pt}
        Assuming you are an expert in the field of mixed-integer programming, please formulate the following problem as a mixed-integer linear programming model and ensure that the model is presented in LaTeX syntax. 

        \hl{Haus Toys can make and sell toy trucks, toy aeroplanes, toy ships, and toy trains. The profit for each truck sold is \$5, each aeroplane sold is \$10, each ship sold is \$8, and each train sold is \$7. How many of each type of toys should Haus Toys make to maximise the profit? 
        
        There are 890 units of timber available. The amount of timber required is 12 units for each truck, 20 units for each aeroplane, 15 units for each ship, and 10 units for each train.
        
        The are 500 units of steel available. The amount of steel required is 3 units for each aeroplane, 5 units for each ship, 4 units for each train, and 6 units for each truck.
        
        If Haus Toys makes trucks, then they will not make trains.
        
        If they make ships, however, they will also make aeroplanes.
        
        The number of toy ships made cannot exceed the number of trains made. 
        
        Haus Toys must also make at least twice as many ships as aeroplanes.}
    \end{mdframed}
\end{table}

\begin{table}[h!]
    \centering
    \caption{An example of the prompt for variable identification}
    \label{fig:prompt for variable identification}
    \begin{mdframed}
        \small
        \setlength{\parskip}{5pt}
        You need to formulate the following given problem as a mixed integer programming (MIP) model. But now you just need to define the decision variables.  
  
        MIP Problem Description: \hl{Haus Toys can make and sell toy trucks, toy aeroplanes, toy ships, and toy trains. The profit for each truck sold is \$5, each aeroplane sold is \$10, each ship sold is \$8, and each train sold is \$7. How many of each type of toys should Haus Toys make to maximise the profit? There are 890 units of timber available. The amount of timber required is 12 units for each truck, 20 units for each aeroplane, 15 units for each ship, and 10 units for each train. There are 500 units of steel available. The amount of steel required is 3 units for each aeroplane, 5 units for each ship, 4 units for each train, and 6 units for each truck. If Haus Toys makes trucks, then they will not make trains. If they make ships, however, they will also make aeroplanes. The number of toy ships made cannot exceed the number of trains made. Haus Toys must also make at least twice as many ships as aeroplanes.}
        
        Please first identify the continuous variables and/or non-binary integer variables in the mixed-integer linear programming problem and provide them as a Python list. If you can find either of these two types of variables, there is no need to find binary variables and include binary variables as part of your answer. If you do not think there are any continuous or integer variables involved in the problem description, then give as an answer a Python list consisting only of the binary variables that you identified as the decision variables for the problem. 
        
        A binary variable usually represents a decision about whether a single activity is carried out or not, such as whether a person does something or not, or whether a class of products is produced or not. Binary variables must be named by starting with ``bi\_" and linking ``bi\_" to the name of the entity it represents. 
        
        An integer variable or continuous variable should describe the decision to be made, usually the quantity to be determined. An integer variable or continuous variable often corresponds to some parameters that represent attributes such as profit, cost, resource consumption per unit, etc. Your answer should only include one list of variable names. Your answer does not need to specify the type of variable and provide any explanation. The variable needs to be named with the name of the quantity it represents and stored in a list. For example, the integer variable representing the quantity of vanilla cakes to be made is named vanilla\_cakes. 
        
        Note that your answer must contain only a list consisting entirely of non-binary integer variables and continuous variables, or a list consisting entirely of binary variables, not a mixed list. 
        
        You need to be aware that the (weighted) sums of some quantities, such as total amount of products, total space, total production time, total amount of resources used, and total cost, are not directly considered as a decision variable if they can be expressed as a (weighted) sum of the variables representing those quantities.
        
    \end{mdframed}

\end{table}

\begin{table}[h!]
    \centering
    \caption{An example of the prompt for objective function generation}
    \label{fig:prompt for obj}
    \begin{mdframed}
        \small
        \setlength{\parskip}{5pt}
        You need to formulate the following given problem as a mixed integer programming (MIP) model. But now you just need to define the objective function. 

        Full problem description to give you context: \hl{Haus Toys can make and sell toy trucks, toy aeroplanes, toy ships, and toy trains. The profit for each truck sold is \$5, each aeroplane sold is \$10, each ship sold is \$8, and each train sold is \$7. How many of each type of toys should Haus Toys make to maximise the profit? There are 890 units of timber available. The amount of timber required is 12 units for each truck, 20 units for each aeroplane, 15 units for each ship, and 10 units for each train. The are 500 units of steel available. The amount of steel required is 3 units for each aeroplane, 5 units for each ship, 4 units for each train, and 6 units for each truck. If Haus Toys makes trucks, then they will not make trains. If they make ships, however, they will also make aeroplanes. The number of toy ships made cannot exceed the number of trains made.  Haus Toys must also make at least twice as many ships as aeroplanes.}
        
        Description of Objective Function: \hl{Haus Toys can make and sell toy trucks, toy aeroplanes, toy ships, and toy trains. The profit for each truck sold is \$5, each aeroplane sold is \$10, each ship sold is \$8, and each train sold is \$7. How many of each type of toys should Haus Toys make to maximise the profit?}  
        
        Please select the relevant variables from the continuous (or integer) variables \hl{[`trucks', `aeroplanes', `ships', `trains']} to build the objective function. Do not change the names of the variables in the list when generating expressions. You must only use the relevant constants from the description of the objective function as the parameters or coefficients of the resulting expression. The parameters or coefficients may involve additional calculations or conversions of units. The resulting expression will not necessarily involve all the variables in the variable list provided. You do not need to model other constraints in the problem.
        
        Your answer needs to indicate whether the objective function should be maximized or minimized. Therefore, your answer should only include the words ``Maximize" or ``Minimize" and an expression.
        
        Please provide the requested information without any extra or unnecessary details and explanations beyond what is explicitly asked.
        
        If the objective function is about profit and the value of the product's profit is given directly, then there is no need to use the difference between the selling price and the cost of the product to calculate the profit value. In this case, the total profit is the sum of the individual profits, independent of costs. If the profit value of a product is not explicitly given, the difference between the selling price and the cost of the product should be used to calculate the profit value.
        
    \end{mdframed}

\end{table}

\begin{table}[h!]
    \centering
    \caption{An example of the prompt for type 3 constraint generation}
    \label{fig:prompt for type3}
    \begin{mdframed}
        \small
        \setlength{\parskip}{5pt}
        You need to formulate the following given problem as a mixed integer programming (MIP) model. But now you just need to define a single constraint of the model.

        Full problem description to give you context: \hl{Haus Toys can make and sell toy trucks, toy aeroplanes, toy ships, and toy trains. The profit for each truck sold is \$5, each aeroplane sold is \$10, each ship sold is \$8, and each train sold is \$7. How many of each type of toys should Haus Toys make to maximise the profit?  There are 890 units of timber available. The amount of timber required is 12 units for each truck, 20 units for each aeroplane, 15 units for each ship, and 10 units for each train. The are 500 units of steel available. The amount of steel required is 3 units for each aeroplane, 5 units for each ship, 4 units for each train, and 6 units for each truck. If Haus Toys makes trucks, then they will not make trains. If they make ships, however, they will also make aeroplanes. The number of toy ships made cannot exceed the number of trains made.  Haus Toys must also make at least twice as many ships as aeroplanes.}
        
        Constraint Description: \hl{There are 890 units of timber available. The amount of timber required is 12 units for each truck, 20 units for each aeroplane, 15 units for each ship, and 10 units for each train.} 
        
        This constraint represents an upper bound on the weighted sum of decision variables. The expression for the constraint has the format ``weighted sum of variables \textless= constant representing the upper bound". Each variable in the constraint inequality will be multiplied by a parameter (weight). This parameter must be a constant mentioned in the constraint description above. Please find the correct parameter in the description that corresponds to each variable to substitute into the inequality based on the variable name. 
        
        Please use the variables \hl{[`trucks', `aeroplanes', `ships', `trains']} to model this constraint without making any alterations to the variable names. You do not need to model other constraints or the objective function in the problem. The resulting expression will not necessarily involve all the variables in the variable list provided. Ensure that the variables in the generated expression completely retain their original names from the list.
           
        Use the symbols ``\textgreater='', ``\textless='', and ``='' to denote greater than or equal to, less than or equal to, and equal to, respectively.
            
        Your answer must only be a mathematical expression and do not provide any explanation.

    \end{mdframed}
\end{table}

\begin{table}[h!]
    \centering
    \caption{An example of the prompt for type 13 logic constraint generation}
    \label{fig:prompt for type13}
    \begin{mdframed}
        \small
        \setlength{\parskip}{5pt}
        You need to formulate the following given problem as a mixed integer programming (MIP) model. But now you just need to define a single logic constraint of the model.

        Full problem description to give you context: \hl{Haus Toys can make and sell toy trucks, toy aeroplanes, toy ships, and toy trains. The profit for each truck sold is \$5, each aeroplane sold is \$10, each ship sold is \$8, and each train sold is \$7. How many of each type of toys should Haus Toys make to maximise the profit?  There are 890 units of timber available. The amount of timber required is 12 units for each truck, 20 units for each aeroplane, 15 units for each ship, and 10 units for each train. The are 500 units of steel available. The amount of steel required is 3 units for each aeroplane, 5 units for each ship, 4 units for each train, and 6 units for each truck. If Haus Toys makes trucks, then they will not make trains. If they make ships, however, they will also make aeroplanes. The number of toy ships made cannot exceed the number of trains made.  Haus Toys must also make at least twice as many ships as aeroplanes.}
        
        Constraint Description: \hl{If Haus Toys makes trucks, then they will not make trains.} 
        
        This constraint belongs to a subtype of logic constraints. Please use the binary variables \hl{[`bi\_trucks', `bi\_aeroplanes', `bi\_ships', `bi\_trains']} to model this logic constraint according to the following guidance. Do not change the names of the variables in the list when generating expressions. You do not need to model other constraints or the objective function in the problem. The resulting expression will not necessarily involve all the variables in the list.
        
        Use the symbols ``\textgreater=", ``\textless=", and ``=" to denote greater than or equal to, less than or equal to, and equal to, respectively.
        
        Template: Consider Statements A and B with truth value represented binary variables a, b respectively - with 1 representing a statement is true and 0 otherwise. 
        In natural language descriptions, this type of constraint often contains a format like ``At most one of A and B" or ``If A then not B" or ``Either A or B or neither (but not both)", which corresponds to the mathematical formula a + b \textless= 1. Please match the variables in the given constraint description with the variables in the template to generate the correct mathematical expression. 
        
        Your answer can only be a mathematical expression and do not provide any explanation.

    \end{mdframed}

\end{table}

\begin{table}[h!]
    \centering
    \caption{An example of the prompt for type 10 logic constraint generation}
    \label{fig:prompt for type10}
    \begin{mdframed}
        \small
        \setlength{\parskip}{5pt}
        You need to formulate the following given problem as a mixed integer programming (MIP) model. But now you just need to define a single logic constraint of the model.

        Full problem description to give you context: 
        \hl{Haus Toys can make and sell toy trucks, toy aeroplanes, toy ships, and toy trains. The profit for each truck sold is \$5, each aeroplane sold is \$10, each ship sold is \$8, and each train sold is \$7. How many of each type of toys should Haus Toys make to maximise the profit?  There are 890 units of timber available. The amount of timber required is 12 units for each truck, 20 units for each aeroplane, 15 units for each ship, and 10 units for each train. The are 500 units of steel available. The amount of steel required is 3 units for each aeroplane, 5 units for each ship, 4 units for each train, and 6 units for each truck. If Haus Toys makes trucks, then they will not make trains. If they make ships, however, they will also make aeroplanes. The number of toy ships made cannot exceed the number of trains made.  Haus Toys must also make at least twice as many ships as aeroplanes.}
        
        Constraint Description: \hl{If they make ships, however, they will also make aeroplanes.} 
        
        This constraint belongs to a subtype of logic constraints. Please use the binary variables \hl{[`bi\_trucks', `bi\_aeroplanes', `bi\_ships', `bi\_trains']} to model this logic constraint according to the following guidance. Do not change the names of the variables in the list when generating expressions. You do not need to model other constraints or the objective function in the problem. The resulting expression will not necessarily involve all the variables in the list.
        
        Use the symbols ``\textgreater=", ``\textless=", and ``=" to denote greater than or equal to, less than or equal to, and equal to, respectively.
        
        Template: Consider Statements A and B with truth value represented binary variables a, b respectively - with 1 representing a statement is true and 0 otherwise. 
        In natural language descriptions, this type of constraint often contains a format like ``If A then B" or ``If not B then not A", which corresponds to the mathematical formula a \textless= b. Please match the variables in the given constraint description with the variables in the template to generate the correct mathematical expression. 
        
        Your answer can only be a mathematical expression and do not provide any explanation.
    \end{mdframed}

\end{table}

\begin{table}[h!]
    \centering
        \caption{An example of the prompt for type 9 constraint generation}
    \label{fig:prompt for type9}
    \begin{mdframed}
        \small
        \setlength{\parskip}{5pt}
        You need to formulate the following given problem as a mixed integer programming (MIP) model. But now you just need to define a single constraint of the model.

        Full problem description to give you context: \hl{Haus Toys can make and sell toy trucks, toy aeroplanes, toy ships, and toy trains. The profit for each truck sold is \$5, each aeroplane sold is \$10, each ship sold is \$8, and each train sold is \$7. How many of each type of toys should Haus Toys make to maximise the profit?  There are 890 units of timber available. The amount of timber required is 12 units for each truck, 20 units for each aeroplane, 15 units for each ship, and 10 units for each train. The are 500 units of steel available. The amount of steel required is 3 units for each aeroplane, 5 units for each ship, 4 units for each train, and 6 units for each truck. If Haus Toys makes trucks, then they will not make trains. If they make ships, however, they will also make aeroplanes. The number of toy ships made cannot exceed the number of trains made.  Haus Toys must also make at least twice as many ships as aeroplanes.}
        
        Constraint Description: \hl{The number of toy ships made cannot exceed the number of trains made.} 
        
        This constraint is a comparison constraint between two variables, expressed in a mathematical formula similar to x \textless= b*y, where x and y are variables and b is a positive constant. 
        
        Please use the variables \hl{[`trucks', `aeroplanes', `ships', `trains']} to model this constraint without making any alterations to the variable names. You do not need to model other constraints or the objective function in the problem. The resulting expression will not necessarily involve all the variables in the variable list provided. Ensure that the variables in the generated expression completely retain their original names from the list.
        
        Use the symbols ``\textgreater=", ``\textless=", and ``=" to denote greater than or equal to, less than or equal to, and equal to, respectively.
        
        Your answer must only be a mathematical expression and do not provide any explanation.
        
    \end{mdframed}
\end{table}

\end{document}